\documentclass[12pt]{amsart}
\usepackage{verbatim}
\usepackage{amsmath,amsfonts,amsthm,amsopn,cite,mathrsfs}
\usepackage{epsfig,verbatim}
\usepackage{subfigure}
\include{macros}

\usepackage{mathdots}

\usepackage[pdftex,
colorlinks, 
bookmarksnumbered, 
bookmarksopen, 
linkcolor=blue, 
citecolor=purple, 
urlcolor=orange, 
]{hyperref}

\usepackage[english, french]{babel}

\usepackage[latin1]{inputenc}

\usepackage{latexsym}
\usepackage{mathabx}
\usepackage{calligra}
\usepackage{graphicx}
\usepackage{subfigure}
\usepackage[justification=centering]{caption}
\usepackage{xcolor}
\usepackage{doi}

\usepackage{url}

\newcommand\bcdot{\ensuremath{%
  \mathchoice%
   {\mskip\thinmuskip\lower0.2ex\hbox{\scalebox{1.5}{$\cdot$}}\mskip\thinmuskip}}%
   {\mskip\thinmuskip\lower0.2ex\hbox{\scalebox{1.5}{$\cdot$}}\mskip\thinmuskip}%
   {\lower0.3ex\hbox{\scalebox{1.2}{$\cdot$}}}%
   {\lower0.3ex\hbox{\scalebox{1.2}{$\cdot$}}}%
}

\newtheorem*{thm*}{Theorem}
\numberwithin{thm}{section}

\theoremstyle{definition}

\renewcommand{\l}{\lambda}							
\newcommand{\nn}{\nonumber}								
\renewcommand{\asymp}{\sim}

\newcommand{\R}{\mathbb{R}}								

\newcommand{\Abs}[1]{\left| #1\right|}						
\newcommand{\bracket}[2]{\left\langle #1 , #2\right\rangle}						

\newcommand{\Lie}[1]{\mathfrak{#1}} 		
\newcommand{\Rspan}[1]{\mathbb{R}\text{-}\mathrm{span}\{ #1\}} 			
\newcommand{\Liez}[1]{\mathfrak{z}(\mathfrak{#1})}						
\newcommand{\coad}{\mathop{\mathrm{ad}^*}}							
\newcommand{\RO}{\mathcal{R}}							

\newcommand{\HS}{\mathcal{H}}			
\newcommand{\RS}{\mathcal{H}_\pi}				
\newcommand{\Pf}{\mathrm{Pf}}						
\newcommand{\rp}{\kappa}										

\renewcommand{\H}{\mathbf{H}_n}					
\newcommand{\h}{\mathfrak{h}_n}						
\newcommand{\HG}[2]{\mathbf{H}_{#1, #2}}						
\newcommand{\HA}[2]{\mathfrak{h}_{#1, #2}}							
\newcommand{\gv}[1]{\tilde{#1}}						

\renewcommand{\L}[2]{L^{#1}(#2)}					
\newcommand{\SF}[1]{\mathscr{S}(\mathbb{R}^{#1})}			
\newcommand{\SFG}[1]{\mathscr{S}(#1)}							

\newcommand{\SL}{\mathcal{L}}	
\newcommand{\HO}{\mathcal{Q}_{\R^n}}		
\newcommand{\HHO}{\mathcal{Q}_{\H}}			

\begin{document}

\centerline{}

\selectlanguage{english}
\title[The Harmonic Oscillator on the Heisenberg Group]{The Harmonic Oscillator on the Heisenberg Group \\ -- \\ L'oscillateur harmonique sur le group de Heisenberg}


\selectlanguage{english}

\author[D. Rottensteiner]{David Rottensteiner}
\address{
David Rottensteiner:
\endgraf Department of Mathematics: Analysis,
Logic and Discrete Mathematics
\endgraf Ghent University
\endgraf Krijgslaan 281, S8
\endgraf 9000 Gent
\endgraf Belgium}
\email{david.rottensteiner@ugent.be}

\author[M. Ruzhansky]{Michael Ruzhansky}
\address{
  Michael Ruzhansky:
  \endgraf
  Department of Mathematics: Analysis,
Logic and Discrete Mathematics
  \endgraf Ghent University
  \endgraf Krijgslaan 281, S8
\endgraf 9000 Gent
\endgraf Belgium
  \endgraf
  and
  \endgraf
  School of Mathematical Sciences
    \endgraf
    Queen Mary University of London
    \endgraf
    Mile End Road 
    \endgraf London E1 4NS
  \endgraf United Kingdom
  \endgraf
  }
\email{ruzhansky@gmail.com}

\subjclass[2010]{35R03, 35P20}
\keywords{Harmonic oscillator, Heisenberg group, Dynin-Folland group, sub-Laplacian, eigenvalue distribution, counting function, eigenfunctions, orthonormal basis}

\begin{abstract}
\selectlanguage{english}

In this note we present a notion of harmonic oscillator on the Heisenberg group $\mathbf{H}_n$ which forms the natural analogue of the harmonic oscillator on $\mathbb{R}^n$ under a few reasonable assumptions: the harmonic oscillator on $\mathbf{H}_n$ should be a negative sum of squares of operators related to the sub-Laplacian on $\mathbf{H}_n$, essentially self-adjoint with purely discrete spectrum, and its eigenvectors should be smooth functions and form an orthonormal basis of $L^2(\mathbf{H}_n)$. This approach leads to a differential operator on $\mathbf{H}_n$ which is determined by the (stratified) Dynin-Folland Lie algebra. We provide an explicit expression for the operator as well as an asymptotic estimate for its eigenvalues.

\vskip 0.5\baselineskip

\selectlanguage{french}
\noindent{\bf R\'esum\'e} \vskip 0.5\baselineskip \noindent
Dans cette note, nous pr\'esentons une notion d'oscillateur harmonique sur le groupe de Heisenberg $\mathbf{H}_n$ qui forme l'analogue naturel de l'oscillateur harmonique sur $\mathbb{R}^n$ sous quelques hypoth\`eses raisonnables: le l'oscillateur harmonique sur $\mathbf{H}_n $ devra\^it \^etre une somme n\'egative de carr\'es d'op\'erateurs li\'ee au sous-laplacien sur $\mathbf{H}_n$, \^etre essentiellement auto-adjoint avec un spectre purement discret, et les vecteurs propres doivent former une base orthonorm\'ee de $L^2(\mathbf{H} _n)$. Cette approche conduit \`a un op\'erateur diff\'erentiel sur $\mathbf{H}_n$ qui est d\'etermin\'e par l'alg\`ebre de Dynin-Folland de Lie (stratifi\'ee). Nous fournissons une expression explicite pour l'op\'erateur ainsi qu'une estimation asymptotique pour ses valeurs propres.

\end{abstract}

\maketitle

\selectlanguage{english}

\vspace{-0.5cm}

\section{Introduction}

The aim of this note is to introduce a canonical harmonic oscillator on the Heisenberg group $\H$. For $n = 1$ and exponential coordinates $(t_1, t_2, t_3) \in \R^3 \cong \mathbf{H}_1$ the harmonic oscillator we propose is explicitly given by
	\begin{align*}
			\mathcal{Q}_{\mathbf{H}_1} = - \bigl( \partial_{t_1}^2 + \partial_{t_2}^2 \bigr) - \frac{1}{4} \bigl ({t_1}^2 + {t_2}^2 \bigr ) \hspace{1pt} \partial_{t_3}^2 + \bigl ( t_1 \hspace{1pt} \partial_{t_2} - t_2 \hspace{1pt} \partial_{t_1} \bigr ) \hspace{1pt} \partial_{t_3} + 4 \pi^2  \hspace{1pt} t_3^2.
	\end{align*}
 Our approach is motivated by the following three realizations of the classical harmonic oscillator $\HO$ on $\R^n$:
\begin{itemize}
	\item[(R1)] the negative sum of squares $- \Delta + 4 \pi^2 \Abs{t}^2$ of partial derivatives of order 1 and coordinate multiplication operators;
	\item[(R2)] the Weyl and Kohn-Nirenberg quantizations on $\R^n$ of the symbol $\sigma(t, \xi) := 4 \pi^2 (\Abs{\xi}^2 + \Abs{t}^2)$ with $t, \xi \in \R^n$;
	\item[(R3)] the image $d\rho_1(-\SL_{\H})$ of the negative sub-Laplacian $-\SL_{\H}$ on $\H$ under the infinitesimal Schr\"{o}dinger representation $d\rho_1$ (of Planck's constant equal to $1$) of the Heisenberg Lie algebra $\h$\footnote{extended to the universal enveloping algebra $\mathfrak{u}(\h)$.}, for $n=1$.
\end{itemize}
\medskip
The operator $\HO$ is usually defined by the expression $- \Delta + 4 \pi^2 \Abs{t}^2$, or some scaled version of it.\footnote{The factor $4 \pi^2$ is due to our choice of realising the Schr\"{o}dinger representation; our expression agrees with the versions in Folland~\cite{FollPhSp} or Stein~\cite{Stein} up to scaling.} However, the Schr\"{o}dinger representation $\rho_1$ of $\H$ acting on $\L{2}{\R^n}$ and the associated Lie algebra representation, naturally acting on $\SF{n}$, clearly relate each of the realizations (R1) - (R3) to the others. Moreover, each of these realizations features a sum of squares. It ought therefore to be natural to assume that similar realizations should be available for the canonical harmonic oscillator on $\H$.

The special role of the Heisenberg Lie algebra $\h$ in this context is not coincidental: it is precisely the Lie algebra which is generated by the partial derivatives $\partial_{t_j}$ and the multiplication operators for the coordinate functions $t_k$, $j,k = 1, \ldots, n$. It is well known that $\h$ is stratified, therefore permits a (canonical) homogeneous structure, and that the sums of squares in the identities above are essentially related to the first stratum of $\h$.

An operator on $\H$ satisfying criteria analogous to (R1) - (R3) should clearly involve left-invariant (or alternatively right-invariant) vector fields on $\H$, which are uniquely determined by some vectors in $\h$, and a scalar potential expressed in terms of the coordinate functions on $\H$. It ought therefore to be natural to study the Lie algebra generated by the standard basis of left-invariant vector fields, here denoted by $X_1, \ldots, X_{2n+1}$, and the multiplication operators defined by the coordinates $t_1, \ldots, t_{2n+1}$ on $\h \cong \R^{2n+1}$ which determine the coordinates in which the vector fields are written.
The resulting Lie algebra and its representation theory were first studied in Dynin~\cite{Dyn1}, and in more detail in Folland~\cite{FollMeta}. This Lie algebra, which we shall call the Dynin-Folland Lie algebra, is in fact stratified and thus admits a sub-Laplacian. Endowed with the canonical homogeneous structure arising from the stratification, the Dynin-Folland Lie algebra together with its associated connected, simply connected Lie group, the group's generic irreducible unitary representations, and the associated negative sub-Laplacian (a positive Rockland operator) give rise to \textit{the harmonic oscillator on the Heisenberg group}. We provide a concrete formula for this operator and describe the asymptotic growth of its eigenvalues, using results by ter Elst and Robinson ~\cite{tERo}.

\section{The Dynin-Folland Group and its Representations}

\subsection{The Dynin-Folland Lie Algebra}

In order to present our results, in this and the next sections we fix the notation and recall the fundamental results about the Dynin-Folland Lie algebra and group and its generic unitary irreducible representations due to Dynin~\cite{Dyn1} and Folland~\cite{FollMeta}. For more details we refer to \cite{FollMeta, Ro14, FiRoRu}.

We choose the usual exponential coordinates for the Heisenberg group $\H$, thus express the group law by
	\begin{align}
		 (t_{2n+1}, t_{2n}, &\ldots, t_1) (t'_{2n+1}, t'_{2n}, \ldots, t'_1) \nn \\ &= \Bigl( t_{2n+1} + t'_{2n+1} + \frac{1}{2} \sum_{j=1}^n (t_j t'_{n+j} - t'_j t_{n+j}), t_{2n} + t'_{2n}, \ldots, t_1 + t'_1 \Bigr).  \nonumber
	\end{align}
We can also group the variables as $\gv{t}_3 := t_{2n+1}, \gv{t}_2 := (t_{2n}, \ldots, t_{n+1}), \gv{t}_1 := (t_n, \ldots, t_1)$ and rewrite the group law as
	\begin{align}
		(\gv{t}_3, \gv{t}_2, \gv{t}_1) (\gv{t}'_3, \gv{t}'_2, \gv{t}'_1) = \Bigl( \gv{t}_3 + \gv{t}'_3 + \frac{1}{2} \bigl( \bracket{\gv{t}_1}{\gv{t}'_2} - \bracket{\gv{t}_2}{\gv{t}'_1} \bigr), \gv{t}_2 + \gv{t}'_2, \gv{t}_1 + \gv{t}'_1 \Bigr). \label{GrLawHgv}
	\end{align}
In these coordinates, one can realise the Schr\"{o}dinger representation $\rho_\rp$ of formal dimension $\Abs{\rp}^n$, $\rp \in \R \setminus \{ 0 \}$, on $f \in \L{2}{\R^n}$ as
	\begin{align}
		\bigl( \rho_{\rp}(\gv{t}_3, \gv{t}_2, \gv{t}_1)f \bigr)(\gv{t}'_1) = e^{2 \pi i \rp \bigl( \gv{t}_3 + \frac{1}{2} \bracket{\gv{t}_1}{\gv{t}_2} + \bracket{\gv{t}_2}{\gv{t}'_1} \bigr)} \hspace{1pt} f(\gv{t}'_1 + \gv{t}_1). \nn
	\end{align}

The real Lie algebra $\widetilde{\Lie{g}}$ of operators generated by the left-invariant vector fields on $\H$
$$X_j = \partial_{t_j} - \frac{1}{2} t_{n+j} \partial_{t_{2n+1}}, \hspace{10pt} X_{n+j} = \partial_{t_{n+j}} + \frac{1}{2} t_j \partial_{t_{2n+1}}, \hspace{10pt} X_{2n+1} = \partial_{t_{2n+1}}, \hspace{5pt} j =1, \ldots, n,$$
and by the multiplication operators $Y_k = 2 \pi i t_k$ for $k = 1, \ldots, 2n+1$ is $3$-step nilpotent and, as a vector space, isomorphic to $\R \times \R^{2n+1} \times \h \cong \R^{4n+3}$. If we denote by $Z$ the multiplication by the constant $2 \pi i$ and identify the operators $Z, Y_1, \ldots, X_{2n+1}$ with the standard basis vectors in $\R^{4n+3}$, the isomorphism is realized by equipping $\R^{4n+3}$ with the Lie bracket defined by
\begin{equation}
	\begin{array}{rclrclrcl}
	[X_j, X_{n+j}] &=& X_{2n+1}, \hspace{10pt} & [X_j, Y_{2n+1}] &=& -\frac{1}{2} Y_{n+j}, \hspace{10pt} & [X_{n+j}, Y_{2n+1}] &=& \frac{1}{2} Y_j, \\
	&& & && & [X_k, Y_k] &=& Z,
	\end{array} \label{LieBr}
\end{equation}
for $j = 1, \ldots, n$ and $k = 1, \ldots, 2n+1$, and vanishing brackets otherwise. We will denote this Lie algebra by $\HA{n}{2}$ and refer to it as the Dynin-Folland Lie algebra. The connected, simply connected Lie group obtained by exponentiation will be referred to as the Dynin-Folland group and denoted by $\HG{n}{2}$. The Lie bracket relations immediately reveal that the Lie sub-algebra generated by $Z, Y_1, \ldots, Y_{2n+1}$ is Abelian, and hence  $\HG{n}{2}$ can be viewed as a semi-direct product of the form $\R^{2n+2} \rtimes \H$.
Using exponential coordinates and identifying any element of $\HG{n}{2}$ with its corresponding coordinate vector $\bigl(z, y_1, \ldots, y_{2n+1}, x_{2n+1}, \ldots, x_1 \bigr) =: (z, y, x) \in \R \times \R^{2n+1} \times \R^{2n+1}$, the $\HG{n}{2}$-group law can be expressed by
	\begin{align}
		(z, y, x) (z', y', x') = \Bigl( z + z' + \frac{1}{2} \bigl( &\bracket{x}{y'} - \bracket{x}{y'} \bigr), \\
		&y + y' + \frac{1}{4} \bigl( \coad(x)y' - \coad(y)x' \bigr), x \cdot x' \Bigr) \nonumber
	\end{align}
provided we denote by $x \cdot x'$ the $x''$-coordinates of $(0, 0, x'') = (0, 0, x) \cdot (0, 0, x') \in \H \leq \HG{n}{2}$ and by $\coad$ the coadjoint representation of $\h \cong \R^{2n+1}$ on $\h^* \cong \R^{2n+1}$ given by $$\coad(\gv{t}_3, \gv{t}_2, \gv{t}_1)(\gv{t}'_3, \gv{t}'_2, \gv{t}'_1) = (0, -\gv{t}_3 \gv{t}'_1, \gv{t}_3 \gv{t}'_2).$$

\subsection{Stratification and Unitary Irreducible Representations}

The Lie bracket \eqref{LieBr} admits the stratification
\begin{align*}
        	\Lie{g}_3 := \R Z, \hspace{10pt} \Lie{g}_2 := \Rspan{Y_1, \ldots, Y_{2n}, X_{2n+1}}, \hspace{10pt} \Lie{g}_1 := \Rspan{Y_{2n+1}, X_{2n}, \ldots, X_1},
\end{align*}
which possesses a canonical family of homogeneous dilations $\{ D_r \}_{r > 0}$ on $\HA{n}{2}$ given by
	\begin{equation}
	\begin{array}{rclcl}
		D_r(Z) &=& r^3 Z, && \\
		D_r(Y_k) &=& r^2 Y_k, \hspace{5pt} D_r(X_{2n+1}) &=& r^2 X_{2n+1}, \\
		D_r(X_k) &=& r X_k, \hspace{5pt} D_r(Y_{2n+1}) &=& r Y_{2n+1},
	\end{array}
	\label{DilationsDF}
	\end{equation}
for $k = 1, \ldots, 2n$. The sub-Laplacian $\SL_{\HG{n}{2}}$ induced by the above stratification is the left-invariant differential operator on $\HG{n}{2}$ corresponding to the sum of squares
	\begin{align*}
		X^2_1 + \ldots +X^2_{2n} + Y^2_{2n+1} \in \Lie{u}(\HA{n}{2}).
	\end{align*}
Moreover, for $l := \rp Z^* \in \Liez{g}^*$ with $\rp \in \R \setminus \{ 0 \}$ the matrix representation of the corresponding symplectic form $B_l = l( [\, . \, , \, . \,])$ immediately reveals that $B_l$ is non-degenerate on $\HA{n}{2} / \Lie{z}(\HA{n}{2}) \times \HA{n}{2} / \Lie{z}(\HA{n}{2})$, i.e., up to Plancherel measure zero all unitary irreducible representations are square-integrable modulo the center $Z(\HG{n}{2})$, and its Pfaffian (characterizing the Plancherel measure) is given by $\Pf(l) = \Abs{\rp}^{2n+1}$. These generic representations of $\HG{n}{2}$, denoted by $\pi_{\rp}, \rp \in \R \setminus \{ 0 \}$, can be induced by the characters $\chi_{\rp} := e^{2 \pi i \bracket{\rp Z^*}{\, . \,}}$ of the normal Abelian subgroup $\R^{2n+2} \leq \HG{n}{2}$: for a fixed $\rp \in \R \setminus \{ 0 \}$ the action of $\pi_{\rp}$ on the representation space $\HS_\rp \cong \L{2}{\H}$ is explicitly given by
	\begin{align}
		\bigl( \pi_{\rp}(z, y, x) f \bigr)(t) = e^{2 \pi i \rp z} \hspace{2pt} e^{2 \pi i \rp \bracket{t \cdot \frac{1}{2} x}{y}} \hspace{2pt} f(t \cdot x) \label{GenGrRep}
	\end{align}
for $f \in \L{2}{\H}$, where $t \cdot \frac{1}{2} x$ and $t \cdot x$ again denote the $\H$-group products of the corresponding coordinate vectors.

\section{The Harmonic Oscillator on $\H$}

The representation $\pi := \pi_{1}$ for $\rp = 1$ defined in \eqref{GenGrRep} was the object of interest in Dynin's account ~\cite{Dyn1} since it served the purpose of introducing a Weyl quantization on $\H$. For our definition of the harmonic oscillator $\HHO$ on $\H$ this representation plays the same crucial role as the Schr\"{o}dinger representation does for $\HO$. For this reason the analog of (R3) yields a canonical definition of $\HHO$. The analogs of (R1) and (R2) will be an immediate consequence of our choice.

\begin{def} \label{DefHHO}
For the basis $\{ Z, Y_1, \ldots, X_{2n+1} \}$ of the Dynin-Folland Lie alebra $\HA{n}{2}$ and the representation $\pi = \pi_1 \in \widehat{\mathbf{H}}_{n, 2}$ realized on the representation space $\L{2}{\H}$, we define the \textbf{harmonic oscillator} on $\H$ to be the positive essentially self-adjoint operator
	\begin{align*}
		\HHO :=& d\pi \bigl( -\SL_{\HG{n}{2}} \bigr) = - d\pi \bigl( X_1\bigr)^2 - \ldots - d\pi \bigl( X_{2n} \bigr)^2 - d\pi \bigl( Y_{2n+1}^2 \bigr),
	\end{align*}
whose natural domain includes the space of smooth vectors $\HS^\infty_{\pi} \cong \SFG{\H}$.
\end{def}
 
\medskip

The essentially self-adjoint differential operator $\HHO$ on $\H$ admits the following three realizations:
\begin{itemize}
	\item[(R1')] the differential operator $- \SL_{\H} + 4 \pi^2 \hspace{1pt} t_{2n+1}^2$;
	\item[(R2')] the Dynin-Weyl quantization on $\H$ of the symbol $\sigma(t, \xi) := 4 \pi^2 (\xi_1^2 + \dots + \xi_{2n}^2 + t_{2n+1}^2)$ with $t, \xi \in \R^{2n+1}$;
	\item[(R3')] the element $d\pi(-\SL_{\H}) \in \Lie{u}(\HA{n}{2})$ for the sub-Laplacian $\SL_{\H}$ on $\H$.
\end{itemize}
Since the Lie algebra isomorphism $\HA{n}{2} \to \widetilde{\Lie{g}}$ defined by \eqref{LieBr} is precisly $d\pi$, we immediately have
\begin{equation} \label{Formula}
\begin{array}{rcl}
\HHO &=& - \SL_{\H} + 4 \pi^2 \hspace{1pt} t_{2n+1}^2 \\
		&=& - \bigl( \partial_{t_1} - \frac{1}{2} t_{n+1} \partial_{t_{2n+1}} \bigr)^2 - \ldots - \bigl( \partial_{t_{2n}} + \frac{1}{2} t_n \partial_{t_{2n+1}}\bigr)^2 + 4 \pi^2 \hspace{1pt} t_{2n+1}^2,
\end{array}
\end{equation}
thus (R1'). As for Dynin's Weyl quanization of $\sigma(t, \xi) := 4 \pi^2 (\xi_1^2 + \dots + \xi_{2n}^2 + t_{2n+1}^2)$, it suffices to recall that every monomial in $\xi_j$, $j = 1, \ldots, 2n+1$, is mapped to the monomial of the left-invariant vector field $X_j$, and that multiplication by any monomial in $t_k$, $k = 1, \ldots, 2n+1$, is mapped to the multiplication operator multiplication operators $Y_k = 2 \pi i t_k$. For more details we refer to \cite{Dyn1} and \cite[\S~5]{Ro14}.

Note that equivalently we could realize $\pi$ as a direct summand of the left regular representation of $\HG{n}{2}$, thereby replacing the left-invariant sub-Laplacian $\SL_{\H}$ in \eqref{Formula} by the right-invariant one. The spectral asymptotics, however, would not change (cf.~\cite{tERo}).

\section{Spectral Properties}

The harmonic oscillator $\HHO$ has purely discrete spectrum in $(0, \infty)$ and we obtain the asymptotic growth rate of its eigenvalues of $\HHO$ by employing a powerful method developed in ter Elst and Robinson~\cite{tERo}, which applies to general graded groups. If $G$ is nilpotent but non-Abelian, then for almost every unitary irreducible representation $\pi \in \widehat{G}$ the representation space $\HS_\pi$ is infinite-dimensional.
Moreover, if $G$ is graded, then a left-invariant differential operator $\RO$ on $G$ is said to be a Rockland operator if for every $\pi \in \widehat{G}$ the operator $d\pi(\RO)$ is injective on the space of smooth vectors $\HS^\infty_\pi \subset \HS_\pi$.

Hulanicki, Jenkins and Ludwig~\cite{HuJeLu} showed that if $\RO$ is positive, then for every $\pi \in \widehat{G}$ the operator $d\pi(\RO)$ has purely discrete spectrum in $(0, \infty)$. A concrete description of the spectrum is due to ter Elst and Robinson~\cite{tERo}, who showed that the number of eigenvalues of $d\pi(\RO)$, counted with multiplicities, asymptotically grows like the volumes of certain subsets of the corresponding coadjoint orbit $\mathcal{O}_\pi$. The subsets in question are determined (up to a multiplicative constant) by a (any) homogeneous quasi-norm on $\Lie{g}^*$. Their estimate also gives an asymptotic value for the magnitude of a given eigenvalue.
\medskip

In the case of $G = \HG{n}{2}$ and $\RO = - \SL_{\HG{n}{2}}$, the realization of $\pi = \pi_1 \in \widehat{\mathbf{H}}_{n, 2}$ in $\HS_\pi = \L{2}{\H}$ given by \eqref{GenGrRep} makes these results readily available for $d\pi(\RO) = \HHO$. The choice of a convenient quasi-norm on $\Lie{h}^*_{n , 2}$ and the fact that the coadjoint orbit $\mathcal{O}_\pi$ is flat facilitate the computation of the volumes in question substantially. One can use this to show:
\medskip

\begin{thm*}
The harmonic oscillator $\HHO$ on the Heisenberg group $\H$ has a purely discrete spectrum $\mathrm{spec}(\HHO) \subset (0, \infty)$. The number of its eigenvalues, counted with multiplicities, which are less or equal to $\l > 0$ is asymptotically given by
	\begin{align*}
		N(\l) \asymp \l^{\frac{6n + 3}{2}},
	\end{align*}
and the magnitude of the eigenvalues is asymptotically equal to
	\begin{align*}
		\l_{s} \asymp s^{ \frac{2}{6n+3}} \hspace{5pt} \mbox{ for } \hspace{5pt} s = 1, 2, \ldots.
	\end{align*}
Moreover, the eigenvectors of $\HHO$ are in $\SFG{\H}$ and form an orthonormal basis of $\L{2}{\H}$.
\end{thm*}
\medskip

The power ${\frac{6n + 3}{2}}$ bears a specific relation to the canonical homogeneous structure of $\h$: the nominator $6n + 3$ is the homogenous dimension of the first two strata of $\HA{n}{2}$, i.e., of the subspace $\Lie{g}_1 \, \oplus \, \Lie{g}_2 \subseteq \HA{n}{2}$, while the denominator $2$ is the homogeneous degree of $- \SL_{\HG{n}{2}}$. For a proof of the spectral asymptotics we refer to the preprint~\cite{RoRu18}, especially the proof of Proposition~6.3.
 The eigenvectors are clearly elements of $\bigcap_{k = 1}^\infty \mathrm{dom}(\overline{\HHO}^k) \subseteq \L{2}{\H}$; since $\bigcap_{k = 1}^\infty \mathrm{dom}(\overline{d\pi(\RO)}^k) = \RS^\infty$ (see~\cite[Prop.~2.1]{tERo}), this set coincides with $\SFG{\H}$. The eigenvectors form an orthonormal basis of $\L{2}{\H}$ because $\HHO$ is essentially self-adjoint and $\lim_{s \to \infty} \l_s = \infty$ (see, e.g., Schm\"{u}dgen~\cite[Prop.~5.12]{Schm}).

\section*{Acknowledgments}

Michael Ruzhansky was supported by the FWO Odysseus 1 grant G.0H94.18N: Analysis and Partial Differential Equations, by the EPSRC grant EP/R003025/1 and by the Leverhulme
Grant RPG-2017-151.

David Rottensteiner was supported by the FWO Odysseus 1 grant G.0H94.18N: Analysis and Partial Differential Equations and the Austrian Science Fund (FWF) project [I~3403].

\bibliographystyle{alphaabbr}
\bibliography{Bib_HAO}

\begin{thebibliography}{FRR18}

\bibitem[Dyn75]{Dyn1}
A.~S. Dynin.
\newblock Pseudodifferential operators on the {H}eisenberg group.
\newblock {\em Dokl. Akad. Nauk SSSR}, 225:1245--1248, 1975.

\bibitem[Fol89]{FollPhSp}
G.~B. Folland.
\newblock {\em Harmonic Analysis in Phase Space}.
\newblock Princeton University Press, 1989.

\bibitem[Fol94]{FollMeta}
G.~B. Folland.
\newblock Meta-{H}eisenberg groups.
\newblock In {\em Fourier analysis ({O}rono, {ME}, 1992)}, volume 157 of {\em
  Lecture Notes in Pure and Appl. Math.}, pages 121--147. Dekker, New York,
  1994.

\bibitem[FRR18]{FiRoRu}
V.~Fischer, D.~Rottensteiner, and M.~Ruzhansky.
\newblock {H}eisenberg-{M}odulation {S}pace on the {C}rossroads of {C}oorbit
  {T}heory and {D}ecomposition {S}pace {T}heory.
\newblock {\em Preprint}, 2018.
\newblock \url{https://arxiv.org/abs/1812.07876}.

\bibitem[HJL85]{HuJeLu}
A.~Hulanicki, J.~W. Jenkins, and J.~Ludwig.
\newblock Minimum eigenvalues for positive, {R}ockland operators.
\newblock {\em Proc. Amer. Math. Soc.}, 94(4):718--720, 1985.

\bibitem[Rot14]{Ro14}
D.~Rottensteiner.
\newblock {\em {T}ime-{F}requency {A}nalysis on the {H}eisenberg {G}roup}.
\newblock PhD thesis, Imperial College London, September 2014.

\bibitem[RR18]{RoRu18}
D.~Rottensteiner and M.~Ruzhansky.
\newblock Harmonic and {A}nharmonic {O}scillators on the {H}eisenberg {G}roup.
\newblock {\em Preprint}, 2018.
\newblock \url{https://arxiv.org/abs/1812.09620}.

\bibitem[Sch12]{Schm}
K.~Schm\"{u}dgen.
\newblock {\em Unbounded self-adjoint operators on {H}ilbert space}, volume 265
  of {\em Graduate Texts in Mathematics}.
\newblock Springer, Dordrecht, 2012.

\bibitem[Ste93]{Stein}
E.~M. Stein.
\newblock {\em Harmonic analysis: real-variable methods, orthogonality, and
  oscillatory integrals}, volume~43 of {\em Princeton Mathematical Series}.
\newblock Princeton University Press, Princeton, NJ, 1993.
\newblock With the assistance of Timothy S. Murphy, Monographs in Harmonic
  Analysis, III.

\bibitem[tER97]{tERo}
A.~F.~M. ter Elst and D.~W. Robinson.
\newblock Spectral estimates for positive {R}ockland operators.
\newblock In {\em Algebraic groups and {L}ie groups}, volume~9 of {\em Austral.
  Math. Soc. Lect. Ser.}, pages 195--213. Cambridge Univ. Press, Cambridge,
  1997.

\end{thebibliography}

\end{document}